\documentclass[a4paper,11pt,reqno]{amsart}
\usepackage[T1]{fontenc}
\usepackage{lmodern,amsmath,amssymb,amsthm,mathtools,booktabs}
\usepackage[margin=1in]{geometry}
\usepackage[colorlinks=true,linkcolor=blue,citecolor=blue,urlcolor=blue]{hyperref}
\usepackage{microtype}
\hypersetup{pdftitle={Cyclotomic factorization and fixed-point insertion in colored permutation groups},pdfauthor={Jiang Zeng}}
\numberwithin{equation}{section}
\newtheorem{theorem}{Theorem}[section]
\newtheorem{proposition}[theorem]{Proposition}
\newtheorem{lemma}[theorem]{Lemma}
\newtheorem{corollary}[theorem]{Corollary}
\theoremstyle{definition}
\newtheorem{example}[theorem]{Example}
\newtheorem{remark}[theorem]{Remark}
\DeclareMathOperator{\des}{des}
\DeclareMathOperator{\odes}{odes}
\DeclareMathOperator{\edes}{edes}
\DeclareMathOperator{\Des}{Des}
\DeclareMathOperator{\maj}{maj}
\DeclareMathOperator{\inv}{inv}
\DeclareMathOperator{\negat}{neg}
\DeclareMathOperator{\col}{col}
\DeclareMathOperator{\fmaj}{fmaj}
\DeclareMathOperator{\Fix}{Fix}
\DeclareMathOperator{\fix}{fix}
\newcommand{\E}{E}
\newcommand{\B}{B}
\newcommand{\calB}{\mathcal B}
\newcommand{\uu}{\mathbf u}
\newcommand{\Sn}{\mathfrak S}
\newcommand{\qbinom}[2]{\genfrac{[}{]}{0pt}{}{#1}{#2}}
\title[Cyclotomic factorization and fixed-point insertion]
{Cyclotomic factorization and fixed-point insertion\newline in colored permutation groups}
\author{Jiang Zeng}
\address{College of Mathematics and Physics, Wenzhou University, Wenzhou, China\\
Institut Camille Jordan, Universit\'e Claude Bernard Lyon 1, CNRS, France}
\email{zeng@math.univ-lyon1.fr}
\date{\today}
\begin{document}
\begin{abstract}
We extend cyclotomic Euler--Mahonian factorization to colored
permutation groups for two inversion statistics: Biagioli--Zeng
length and the inversion number of the underlying unsigned
permutation. At a primitive $m$th root and rank $n=mk+i$, each
normalized enumerator is a Hadamard product of a rank-$i$ colored
enumerator and a weighted rank-$k$ type-B enumerator. Separate
color variables are retained. The type-B factors agree for the
two statistics, yielding equidistribution modulo the rank jointly
with descent number, major index and negative-entry number.
We also determine their exact changes under positive fixed-point
insertion and obtain rational kernels retaining descent and major
index. These kernels mark arbitrary fixed-point numbers and give
derangement expansions at every inversion weight. Signed
specializations recover type-A and type-B identities and give
colored flag-major refinements. Inverse statistics and color-sum
subgroups admit related insertion formulas.
\end{abstract}
\maketitle

\section{Introduction and main results}

The joint distribution of descent number, major index and inversion
number is encoded by
\[
 A_n(t,q,p)=\sum_{\sigma\in\Sn_n}
 t^{\des(\sigma)}q^{\maj(\sigma)}p^{\inv(\sigma)}.
\]
Specializing the inversion variable at a root of unity relates signed
and cyclotomic enumeration to ordinary Euler--Mahonian enumeration.
Wachs's involution treats the even signed case \cite{Wachs}, and
Adin, Gessel and Roichman obtained related Mahonian evaluations
\cite{AGR}. Han \cite{Han} proved that, for a primitive $m$th root
$\zeta$ and $n=mk+i$, $0\le i<m$,
\begin{equation}\label{eq:Han}
 \frac{A_{mk+i}(t,q,\zeta)}{(t;q)_{mk+i+1}}
 =\frac{A_i(t,q,\zeta)}{(t;q)_{i+1}}
 \otimes\frac{A_k(t,q^m,1)}{(t;q^m)_{k+1}}.
\end{equation}
Here $(t;q)_N=\prod_{j=0}^{N-1}(1-tq^j)$, with $(t;q)_0=1$,
and $\otimes$ is coefficientwise multiplication of power series in
$t$. In Question 21 of \cite{Han}, Han asks for an extension to
colored permutation groups.

We give two such factorizations, for colored length and for the
inversion number of the underlying unsigned permutation. Their
proofs use finite Gaussian-multinomial formulas and the same
no-carry evaluation at roots of unity. We then study fixed-point
insertion. Its exact inversion increments give rational kernels,
fixed-point-number refinements and derangement formulas at
arbitrary values of the inversion variable.

Write $G(r,n)=G(r,1,n)=C_r\wr\Sn_n$. A colored window is
$\gamma=[\sigma_1^{c_1},\ldots,\sigma_n^{c_n}]$, with
$\sigma\in\Sn_n$ and $0\le c_j<r$; write $|\gamma|=\sigma$.
Our principal descent and major-index statistics use the
Biagioli--Zeng order, specified in Section~\ref{sec:stats}.
For independent color weights $\uu=(u_1,\ldots,u_{r-1})$ and
$u_0=1$, put $w(\gamma)=\prod_j u_{c_j}$ and define
\begin{align}
 \E_n^{(r)}(t,q,p;\uu)
 &=\sum_{\gamma\in G(r,n)}
 t^{\des_G(\gamma)}q^{\maj(\gamma)}p^{\ell_G(\gamma)}w(\gamma),
 \label{eq:Edef}\\
 \mathcal U_n^{(r)}(t,q,p;\uu)
 &=\sum_{\gamma\in G(r,n)}
 t^{\des_G(\gamma)}q^{\maj(\gamma)}p^{\inv(|\gamma|)}w(\gamma).
 \label{eq:Udef}
\end{align}
Both reduce to $A_n$ when $r=1$. For $r=2$ we use
\begin{align}
 \B_n(t,q,p;z)&=\sum_{\beta\in B_n}
 t^{\des_B(\beta)}q^{\maj(\beta)}p^{\ell_B(\beta)}z^{\negat(\beta)},
 \label{eq:Bdef}\\
 I_n(t,q,p;z)&=\sum_{\beta\in B_n}
 t^{\des_B(\beta)}q^{\maj(\beta)}p^{\inv(|\beta|)}z^{\negat(\beta)},
 \label{eq:Idef}\\
 \calB_n(t,q;z)&=\B_n(t,q,1;z)=I_n(t,q,1;z).
 \label{eq:calBdef}
\end{align}
All rank-zero enumerators are one.

\begin{theorem}[Two colored cyclotomic factorizations]\label{thm:main}
Let $r\ge1$, $m\ge2$, $n=mk+i$ with $0\le i<m$, and let
$\zeta$ be a primitive $m$th root of unity. For
$\mathcal F=\E$ or $\mathcal U$, set
\begin{equation}\label{eq:effective}
 \alpha_{\E}=\sum_{c=1}^{r-1}u_c\zeta^{c-1},\qquad
 \alpha_{\mathcal U}=\sum_{c=1}^{r-1}u_c,\qquad
 b_{\mathcal F}=-(-\alpha_{\mathcal F})^m.
\end{equation}
Then
\begin{equation}\label{eq:main}
 \frac{\mathcal F_{mk+i}^{(r)}(t,q,\zeta;\uu)}{(t;q)_{mk+i+1}}
 =\frac{\mathcal F_i^{(r)}(t,q,\zeta;\uu)}{(t;q)_{i+1}}
 \otimes\frac{\calB_k(t,q^m;b_{\mathcal F})}{(t;q^m)_{k+1}}.
\end{equation}
Equivalently, for $\ell\ge0$,
\begin{align}
 [t^\ell]\frac{\mathcal F_{mk+i}^{(r)}(t,q,\zeta;\uu)}{(t;q)_{mk+i+1}}
 &=\left([t^\ell]\frac{\mathcal F_i^{(r)}(t,q,\zeta;\uu)}{(t;q)_{i+1}}\right)
 \notag\\[-2pt]
 &\quad\times
 \bigl([\ell+1]_{q^m}+b_{\mathcal F}[\ell]_{q^m}\bigr)^k.
 \label{eq:maincoeff}
\end{align}
\end{theorem}

Here $[s]_Q=1+Q+\cdots+Q^{s-1}$, with $[0]_Q=0$.
These are formal power-series identities with polynomial
coefficients over $\mathbb Q(\zeta)$. For the length enumerator
we also write $\alpha=\alpha_{\E}$ and $b=b_{\E}$.
The scalar color weights $u_c=a^c$ give
\begin{equation}\label{eq:scalar}
 \alpha_{\E}=a[r-1]_{a\zeta},\qquad
 \alpha_{\mathcal U}=a[r-1]_a.
\end{equation}
When $r=1$, both effective weights vanish and
$\calB_k(t,q;0)=A_k(t,q,1)$, recovering~\eqref{eq:Han}.
For type B the two effective weights agree. This yields the
comparison congruence in Corollary~\ref{cor:rankcongruence}:
\[
 I_n(t,q,p;z)\equiv\B_n(t,q,p;z)\pmod{p^n-1}.
\]

Biagioli--Zeng's extension of the Garsia--Gessel method
\cite{BZ,GG} supplies the finite length coefficient formula.
For underlying inversions, a signed $P$-partition expansion has
decreasing unsigned labels inside each negative block. Their
contribution $p^{\binom b2}$ produces the same type-B bulk
factor after root specialization. Eu, Lin and Lo \cite{ELL}
studied signed Euler--Mahonian identities on colored groups.
The order-two cases of our two type-B formulas recover their
Theorems 5.1 and 5.2. For $m>2$ the inversion weights considered
here are generally not group characters.

Our second construction concerns positive fixed points, meaning
$\sigma_j=j$ and $c_j=0$. Let $\operatorname{dp}(\pi)$ be the
colored derangement obtained by deleting these entries and
standardizing. Write $X_L(\gamma)=\ell_G(\gamma)$ and
$X_I(\gamma)=\inv(|\gamma|)$. The insertion of $g_h$ fixed points
in gap $h$ of a core $\delta$ has the exact increment
\[
 X(I_g(\delta))-X(\delta)=2\sum_hg_h H_h^X(\delta),
 \qquad H_h^L=\rho_h,\quad H_h^I=\eta_h,
\]
where $\rho_h,\eta_h$ are the cut heights in
Section~\ref{sec:insertion}. Theorem~\ref{thm:kernel} gives a
rational kernel retaining descent number, major index and either
inversion statistic. It also applies to arbitrary colored cores;
Theorem~\ref{thm:markedkernel} then marks every positive fixed
point and yields an inclusion-exclusion formula for derangements.

The insertion coordinates are those used by Faliharimalala and
Zeng \cite[Section 3]{FZ}. The ordinary fixed-point height weight
also occurs in Shin--Zeng's Laguerre-history continued fraction
\cite[Corollary 11]{SZ}. We retain the actual descent and major
index on each fiber together with the exact colored inversion
increment. The inverse-statistic refinement is given in
Appendix~\ref{app:inverse}.

Ji and Zhang \cite{JZ} obtain signed Mahonian derangement formulas
using deletion and shuffle identities. Their color-block flag-major
convention is distinguished from the Biagioli--Zeng convention in
Section~\ref{sec:stats}. The colored major-index insertion identity
of Faliharimalala--Zeng and Chen--Gao--He \cite{FZ,CGH}, together
with parity preservation, gives the signed exact-fixed-point
formulas of Section~\ref{sec:derangements}. The same section
compares the natural-order signed products with Eu--Lin--Lo and
Pan--Xu--Zeng \cite{ELL,PXZ}. At higher roots, the insertion
heights are essential: a direct derangement substitution in the
full-group Hadamard factorization already fails in rank two.

Section~\ref{sec:proof} proves the two factorizations and their
specializations. Section~\ref{sec:insertion} develops insertion
and fixed-point marking, and Section~\ref{sec:derangements}
contains the derangement and signed applications.
Section~\ref{sec:restrictions} gives color restrictions and further
questions. The longer inverse and Galois-twisted refinements are
placed in appendices, followed by exact computational checks.

\section{Statistics and coefficient formulas}\label{sec:stats}

\subsection{The Biagioli--Zeng conventions}
We put the following total order on colored letters, including $0$:
\begin{equation}\label{eq:order}
 n^{r-1}<\cdots<n^1<\cdots<1^{r-1}<\cdots<1^1
 <0<1<\cdots<n.
\end{equation}
Color zero is omitted from the notation. For
$\gamma=[\sigma_1^{c_1},\ldots,\sigma_n^{c_n}]$, set
\begin{align}
 \Des_G(\gamma)&=\{j\in\{0,\ldots,n-1\}:\gamma(j)>\gamma(j+1)\},
 \qquad\gamma(0)=0,\label{eq:des}\\
 \des_G(\gamma)&=|\Des_G(\gamma)|,
 \qquad \maj(\gamma)=\sum_{j\in\Des_G(\gamma)}j,\label{eq:maj}\\
 \inv(\gamma)&=\#\{(h,j):1\le h<j\le n,\ \gamma(h)>\gamma(j)\},
 \label{eq:inv}\\
 \ell_G(\gamma)&=\inv(\gamma)+\sum_{j:c_j>0}(\sigma_j+c_j-1),
 \qquad \col(\gamma)=\sum_{j=1}^n c_j.\label{eq:length}
\end{align}
These are the conventions of \cite[Section 2]{BZ}. For $r=2$ they
give the natural order $-n<\cdots<-1<0<1<\cdots<n$ and the usual
type-B descent set and Coxeter length. An unsubscripted $\maj$ refers to this descent set; a descent at
zero contributes to $\des_G$ but contributes zero to $\maj$.
The color-block major index introduced below always carries the
subscript $\prec$.

The group inverse is specified by
\[
 \gamma^{-1}(\sigma_j)=j^{\bar c_j},\qquad
 \bar c_j\equiv-c_j\pmod r,\quad 0\le\bar c_j<r.
\]
All inverse statistics below refer to this group inverse.
We call $j$ a positive fixed point if $\sigma_j=j$ and $c_j=0$,
and put
\[
 \Fix_0(\gamma)=\{j:\sigma_j=j,\ c_j=0\},\qquad
 \fix_0(\gamma)=|\Fix_0(\gamma)|,\qquad
 \mathcal D_n^{(r)}=\{\gamma:\Fix_0(\gamma)=\varnothing\}.
\]
Deleting these positions leaves the same set of absolute values
as positions. Increasing standardization of both, with colors
unchanged, gives $\operatorname{dp}(\gamma)\in\mathcal D_k^{(r)}$.
For signed windows we also write $\Fix_+(\beta)=\{j:\beta_j=j\}$.

\begin{lemma}[Color reduction]\label{lem:fold}
Set $\alpha(p)=\sum_{c=1}^{r-1}u_c p^{c-1}$. Then
\begin{equation}\label{eq:fold}
 \E_n^{(r)}(t,q,p;\uu)=\B_n(t,q,p;\alpha(p)).
\end{equation}
The unsigned-inversion counterpart is
\begin{equation}\label{eq:foldI}
 \mathcal U_n^{(r)}(t,q,p;\uu)
 =I_n(t,q,p;\alpha_0),\qquad \alpha_0=\sum_{c=1}^{r-1}u_c.
\end{equation}
Both reductions preserve the complete positive fixed-point set.
\end{lemma}
\begin{proof}
Replace each uncolored letter $\sigma_j$ by $+\sigma_j$ and each
letter of nonzero color by $-\sigma_j$, obtaining $\beta\in B_n$.
Distinct letters of a colored permutation have distinct absolute
values. Thus their comparisons in~\eqref{eq:order} depend only on
their absolute values and on whether their colors are zero. This
replacement preserves the inversion number and the descent set.
Moreover,
\[
 \ell_G(\gamma)=\ell_B(\beta)+\sum_{j:c_j>0}(c_j-1).
\]
For fixed $\beta$, each negative entry can independently receive any
color $c\in\{1,\ldots,r-1\}$, with weight $u_cp^{c-1}$. Summing
over the fiber above $\beta$ gives $\alpha(p)^{\negat(\beta)}$.
For $r=1$, only windows with no negative entries survive. This proves
\eqref{eq:fold} in all cases. The unsigned permutation and the
positive fixed-point set are unchanged by folding. For the
unsigned-inversion weight the color contribution is simply
$u_c$, which proves~\eqref{eq:foldI}.
\end{proof}

For $u_c=a^c$, the substitution $\alpha(p)=a[r-1]_{ap}$ is the
reduction underlying \cite[Remark 5.2]{BZ}. Lemma~\ref{lem:fold}
records its version with separate variables for the colors.
The restriction to the order~\eqref{eq:order} is essential to this
argument.

\subsection{The two inversion statistics and the color-block order}
For colored windows, $X_L=\ell_G$ and $X_I=\inv(|\gamma|)$ are
different statistics. In type B they satisfy
\begin{equation}\label{eq:lengthparity}
 \ell_B(\beta)\equiv\inv(|\beta|)+\negat(\beta)\pmod2.
\end{equation}
Indeed each Coxeter generator changes both sides modulo two.
By the color reduction,
\begin{equation}\label{eq:coloredparity}
 \ell_G(\gamma)\equiv\inv(|\gamma|)+\col(\gamma)\pmod2.
\end{equation}
This congruence is an identity of statistics for every $r$; when
$r$ is odd, its color-parity factor is not a group character.

We also use the color-block order
\begin{equation}\label{eq:blockorder}
 1^{r-1}\prec\cdots\prec n^{r-1}\prec\cdots
 \prec1^1\prec\cdots\prec n^1\prec0\prec1\prec\cdots\prec n.
\end{equation}
Let $\Des_\prec$, $\des_\prec$ and $\maj_\prec$ be computed
in this order, with a possible descent at zero, and put
\begin{equation}\label{eq:blockflag}
 \fmaj_\prec(\gamma)=r\maj_\prec(\gamma)+\col(\gamma).
\end{equation}
The two major-index conventions used in this paper are:
\begin{center}
\begin{tabular}{lll}
\toprule
Order & Major index & Flag-major index\\
\midrule
Biagioli--Zeng, \eqref{eq:order} & $\maj$ & $r\maj+\col$\\
Color-block, \eqref{eq:blockorder} & $\maj_\prec$ & $\fmaj_\prec$\\
\bottomrule
\end{tabular}
\end{center}
The principal factorization theorem uses the first row. The
Ji--Zhang comparison and the signed flag-major fixed-point
formulas use the second row. The length statistic remains
\eqref{eq:length} in both cases.

\subsection{A coefficient formula for length}
Write $[n]_p!=\prod_{j=1}^n[j]_p$ and
\[
 \qbinom{n}{n_0,\ldots,n_\ell}_p
 =\frac{[n]_p!}{[n_0]_p!\cdots[n_\ell]_p!}
 \quad(n_0+\cdots+n_\ell=n).
\]
This is a polynomial in $p$. Define
\[
 H_n(z,p)=\prod_{v=1}^n(1+zp^v),\qquad H_0(z,p)=1.
\]
The $r=2$ case of the Biagioli--Zeng identity
\cite[Theorem 5.1]{BZ} is
\begin{equation}\label{eq:BZ}
 \sum_{n\ge0}\frac{\B_n(t,q,p;z)x^n}
 {(t;q)_{n+1}[n]_p!H_n(z,p)}
 =\sum_{\ell\ge0}t^\ell
 \left(\prod_{j=0}^{\ell-1}e_p(q^jx)\right)
 \widehat e_{z,p}(q^\ell x),
\end{equation}
where
\[
 e_p(x)=\sum_{n\ge0}\frac{x^n}{[n]_p!},\qquad
 \widehat e_{z,p}(x)=\sum_{n\ge0}
 \frac{x^n}{[n]_p!H_n(z,p)}.
\]
We use~\eqref{eq:BZ} only with indeterminate parameters, before
passing to a polynomial formula.

\begin{proposition}\label{prop:coeff}
For every $n,\ell\ge0$,
\begin{equation}\label{eq:finite}
 [t^\ell]\frac{\B_n(t,q,p;z)}{(t;q)_{n+1}}
 =\sum_{\substack{n_0+\cdots+n_\ell=n\\n_j\ge0}}
 q^{\sum_{j=0}^{\ell}jn_j}
 \qbinom{n}{n_0,\ldots,n_\ell}_p
 \prod_{v=n_\ell+1}^{n}(1+zp^v).
\end{equation}
\end{proposition}
\begin{proof}
Extract $t^\ell x^n$ in~\eqref{eq:BZ}, then multiply by
$[n]_p!H_n(z,p)$. The coefficient from the last exponential has
denominator $[n_\ell]_p!H_{n_\ell}(z,p)$, whereas the other
exponentials have denominators $[n_j]_p!$. The factorials produce
the multinomial coefficient, and the remaining factor is the product
in~\eqref{eq:finite}. Both sides are polynomials in $p,q,z$.
\end{proof}

\begin{remark}\label{rem:zeros}
The quotient notation $H_n/H_{n_\ell}$ is harmless before
specialization, but should not be evaluated as a quotient at a root
of unity: numerator and denominator may both vanish. The finite
product in~\eqref{eq:finite} is the polynomial obtained after
cancellation. All subsequent specializations use this product.
\end{remark}

Setting $p=1$ in~\eqref{eq:finite} and applying the multinomial
theorem gives the weighted type-B Carlitz identity
\begin{equation}\label{eq:CarlitzB}
 \frac{\calB_n(t,q;z)}{(t;q)_{n+1}}
 =\sum_{\ell\ge0}
 \bigl(q^\ell+(1+z)[\ell]_q\bigr)^n t^\ell
 =\sum_{\ell\ge0}
 \bigl([\ell+1]_q+z[\ell]_q\bigr)^n t^\ell.
\end{equation}
In particular, this identity is valid at $z=-1$.

\subsection{A coefficient formula for underlying inversions}

\begin{proposition}\label{prop:coeffI}
For $n,\ell\ge0$,
\begin{align}
 &[t^\ell]\frac{I_n(t,q,p;z)}{(t;q)_{n+1}}\notag\\
 &\quad=\sum_{\substack{a_0+\cdots+a_\ell+b_0+\cdots+b_{\ell-1}=n\\
                         a_j,b_j\ge0}}
 q^{\sum_{j=0}^\ell ja_j+\sum_{j=0}^{\ell-1}jb_j}
 z^{\sum_{j=0}^{\ell-1}b_j}
 p^{\sum_{j=0}^{\ell-1}\binom{b_j}{2}}
 \qbinom{n}{a_0,\ldots,a_\ell,b_0,\ldots,b_{\ell-1}}_p.
 \label{eq:finiteI}
\end{align}
For $\ell=0$ the list of $b$'s is empty and the formula is $1$.
\end{proposition}

\begin{proof}
For a fixed $\beta$, consider sequences
\[
 \ell\ge v_1\ge\cdots\ge v_n\ge0,
\]
with $v_j>v_{j+1}$ whenever $j\in\Des_B(\beta)$ and $j\ge1$,
and with $\ell>v_1$ when $0\in\Des_B(\beta)$.
Their generating function, including the upper bound, is
\begin{equation}\label{eq:compatible}
 \sum_{\ell\ge0}t^\ell\sum_{v\text{ compatible}}
 q^{v_1+\cdots+v_n}
 =\frac{t^{\des_B(\beta)}q^{\maj_B(\beta)}}{(t;q)_{n+1}}.
\end{equation}
Indeed, subtracting the mandatory strict steps leaves $n+1$
independent nonnegative gaps. Their respective weights are
$t,tq,\ldots,tq^n$, and the mandatory steps give the numerator.

At level $j$, the letters of $\beta$ with $v_a=j$ must be
increasing in the signed order. Thus their negative absolute values
are decreasing, followed by increasing positive values. Let $a_j$
and $b_j$ count these positive and negative letters. The top level
$\ell$ cannot contain a negative letter, since its first letter
would then be negative and would violate $\ell>v_1$.

Choose the disjoint sets of absolute labels for the $2\ell+1$
blocks. The generating polynomial of inversions between blocks is
the Gaussian multinomial in~\eqref{eq:finiteI}; its value depends
only on their sizes. Inside a positive block there are no unsigned
inversions, and inside a negative block there are exactly
$\binom{b_j}{2}$. The remaining weights are
$q^{\sum j(a_j+b_j)}z^{\sum b_j}$. Sum over all sizes and use
\eqref{eq:compatible} to obtain the result.
\end{proof}

\section{Cyclotomic factorizations and congruences}\label{sec:proof}

Fix a primitive $m$th root $\zeta$ and first work with the type-B weight $z$.

\begin{lemma}[Multinomial root-of-unity evaluation]\label{lem:Lucas}
Let $\zeta$ be a primitive $m$th root, and let $n=mk+i$ with
$0\le i<m$. Write $n_j=mk_j+i_j$, $0\le i_j<m$, and assume
$\sum_{j=0}^{\ell}n_j=n$. The multinomial
$\qbinom{n}{n_0,\ldots,n_\ell}_{\zeta}$ vanishes unless
$\sum_ji_j=i$. When this condition holds,
\begin{equation}\label{eq:Lucas}
 \qbinom{n}{n_0,\ldots,n_\ell}_{\zeta}
 =\binom{k}{k_0,\ldots,k_\ell}
 \qbinom{i}{i_0,\ldots,i_\ell}_{\zeta}.
\end{equation}
\end{lemma}
\begin{proof}
This is the multinomial form of the root-of-unity evaluation often
called the $q$-Lucas theorem; see \cite{Olive} and
\cite[Theorem 18]{Han}. The vanishing assertion can also be seen
directly. The order of vanishing of $[N]_p!$ at a primitive $m$th
root is $\lfloor N/m\rfloor$. Thus the order for the multinomial is
$k-\sum_jk_j$. Since $\sum_ji_j\equiv i\pmod m$, a carry means
$\sum_ji_j=i+md$ with $d\ge1$, and then $k-\sum_jk_j=d>0$.
For completeness, the no-carry evaluation can be read from the same
factorials. For $h\ge1$,
\[
 \lim_{p\to\zeta}\frac{[hm]_p}{[m]_p}=h.
\]
The factors with indices divisible by $m$ therefore supply
$k!/\prod_jk_j!$. The nonmultiples of $m$ in every complete block
have the same product $\prod_{v=1}^{m-1}[v]_\zeta$; those complete
blocks cancel because $\sum_jk_j=k$. The remaining factors give
$[i]_\zeta!/\prod_j[i_j]_\zeta!$, proving~\eqref{eq:Lucas}.
\end{proof}

\begin{lemma}\label{lem:cycle}
Put $C=1-(-z)^m$. Under the no-carry conditions in
Lemma~\ref{lem:Lucas},
\begin{equation}\label{eq:cycle}
 \prod_{v=n_\ell+1}^{n}(1+z\zeta^v)
 =C^{k-k_\ell}\prod_{v=i_\ell+1}^{i}(1+z\zeta^v).
\end{equation}
\end{lemma}
\begin{proof}
The polynomial identity
\[
 \prod_{v=1}^{m}(1+z\zeta^v)=1-(-z)^m
\]
follows from $\prod_{v=1}^{m}(X-\zeta^v)=X^m-1$.
Since $i_\ell\le i$, the interval of indices from $n_\ell+1$ to
$n$ comprises $k-k_\ell$ complete periods and a final segment
with residues $i_\ell+1,\ldots,i$. Multiplying their contributions
proves~\eqref{eq:cycle} as a polynomial identity in $z$.
\end{proof}

\begin{lemma}[The phase of a negative block]\label{lem:blockphase}
For $v\ge0$ and $0\le s<m$,
\begin{equation}\label{eq:binomialphase}
 \zeta^{\binom{mv+s}{2}}
 =(-1)^{(m-1)v}\zeta^{\binom{s}{2}}.
\end{equation}
\end{lemma}
\begin{proof}
Subtracting $\binom s2$ from the exponent leaves
$\binom{mv}{2}+mvs$. For odd $m$ its phase is one; for even
$m$, use $\zeta^{m/2}=-1$ to obtain $(-1)^v$.
\end{proof}

\begin{theorem}[Two weighted type-B factorizations]\label{thm:B}
For $n=mk+i$, $0\le i<m$, and a primitive $m$th root $\zeta$,
put $b=-(-z)^m$. For $F=\B$ or $I$,
\begin{equation}\label{eq:Bmain}
 \frac{F_{mk+i}(t,q,\zeta;z)}{(t;q)_{mk+i+1}}
 =\frac{F_i(t,q,\zeta;z)}{(t;q)_{i+1}}
 \otimes\frac{\calB_k(t,q^m;b)}{(t;q^m)_{k+1}}.
\end{equation}
\end{theorem}
\begin{proof}
First take $F=\B$.
For fixed $\ell$, evaluate~\eqref{eq:finite} at $p=\zeta$.
By Lemma~\ref{lem:Lucas}, only compositions with
$\sum_jk_j=k$ and $\sum_ji_j=i$ contribute. Conversely, any two
such compositions give a unique composition $n_j=mk_j+i_j$ of $n$.
Use Lemma~\ref{lem:cycle} and write $C=1+b$. The sum factors as
\begin{align}
 & [t^\ell]\frac{\B_{mk+i}(t,q,\zeta;z)}{(t;q)_{mk+i+1}}
 \notag\\
 &=\left(\sum_{\substack{i_0+\cdots+i_\ell=i\\i_j\ge0}}
 q^{\sum_jji_j}\qbinom{i}{i_0,\ldots,i_\ell}_{\zeta}
 \prod_{v=i_\ell+1}^{i}(1+z\zeta^v)\right)\notag\\
 &\hspace{1em}\times
 \left(\sum_{\substack{k_0+\cdots+k_\ell=k\\k_j\ge0}}
 \binom{k}{k_0,\ldots,k_\ell}
 q^{m\sum_jjk_j}C^{k-k_\ell}\right).
 \label{eq:separate}
\end{align}
The first parenthesis is the rank-$i$ coefficient in
Proposition~\ref{prop:coeff}; the bounds $i_j<m$ are automatic
because their sum is $i<m$. The second is, by the ordinary
multinomial theorem,
\[
 \left(C\sum_{j=0}^{\ell-1}q^{mj}+q^{m\ell}\right)^k
 =\bigl([\ell+1]_{q^m}+b[\ell]_{q^m}\bigr)^k.
\]
Formula~\eqref{eq:CarlitzB} identifies this as the $t^\ell$
coefficient of $\calB_k(t,q^m;b)/(t;q^m)_{k+1}$.
This proves~\eqref{eq:Bmain} coefficient by coefficient.
If $C=0$, the same argument applies by polynomial specialization;
zero exponents in the displayed sums contribute the constant $1$.

Now take $F=I$.
In~\eqref{eq:finiteI}, write
\[
 a_j=mA_j+\alpha_j,\qquad b_j=mB_j+\beta_j,
 \qquad 0\le\alpha_j,\beta_j<m.
\]
By Lemmas~\ref{lem:Lucas} and~\ref{lem:blockphase}, surviving terms satisfy
$\sum\alpha_j+\sum\beta_j=i$ and $\sum A_j+\sum B_j=k$.
The small-residue contribution is precisely
$[t^\ell]I_i(t,q,\zeta;z)/(t;q)_{i+1}$, by
Proposition~\ref{prop:coeffI}.

A positive bulk block at level $j$ has weight $q^{mj}$ per unit.
A negative bulk block has weight
$(-1)^{m-1}z^mq^{mj}=bq^{mj}$ per unit, by
\eqref{eq:binomialphase}. Summing the ordinary multinomial over
the $A_j,B_j$ gives
\[
 \left(\sum_{j=0}^{\ell}q^{mj}
       +b\sum_{j=0}^{\ell-1}q^{mj}\right)^k,
\]
which is the coefficient version of~\eqref{eq:Bmain}.
The weighted Carlitz identity~\eqref{eq:CarlitzB} identifies the
bulk factor with the second Hadamard factor in~\eqref{eq:Bmain}.
\end{proof}

\begin{proof}[Proof of Theorem~\ref{thm:main}]
For $\mathcal F=\E$, apply~\eqref{eq:fold} with
$z=\alpha_{\E}$ in Theorem~\ref{thm:B}. For
$\mathcal F=\mathcal U$, apply~\eqref{eq:foldI} with
$z=\alpha_{\mathcal U}$. The same color reduction identifies
the small-rank factor in each case. Formula~\eqref{eq:maincoeff}
follows from~\eqref{eq:CarlitzB}.
\end{proof}

\subsection{Comparison of the two inversion statistics}
\begin{corollary}[Equidistribution modulo the rank]\label{cor:rankcongruence}
For every $n\ge1$,
\begin{equation}\label{eq:rankcongruence}
 I_n(t,q,p;z)\equiv\B_n(t,q,p;z)\pmod{p^n-1}.
\end{equation}
Thus Coxeter length and underlying unsigned inversion have the
same distribution modulo $n$, jointly with
$(\des_B,\maj,\negat)$ on the full group $B_n$.
\end{corollary}
\begin{proof}
If $m\ge2$ divides $n$, the remainder in Theorem~\ref{thm:B}
is zero. Since $I_0=\B_0=1$, their evaluations agree at every
primitive $m$th root. Equality at $p=1$ follows from the
definitions. Hence their difference vanishes at every $n$th root
of unity and is divisible by the monic square-free polynomial
$p^n-1$. Reducing the exponents of $p$ modulo $n$ proves the
coefficient interpretation.
\end{proof}

\begin{corollary}\label{cor:BIdegenerate}
If $(-z)^m=1$, then, for $F=\B$ or $I$,
\begin{equation}\label{eq:BIdegenerate}
 F_{mk+i}(t,q,\zeta;z)=(t;q)_{mk}F_i(tq^{mk},q,\zeta;z).
\end{equation}
In particular this holds for every even $m$ at $z=1$.
\end{corollary}
\begin{proof}
The Carlitz coefficient for the bulk factor is
$([\ell+1]_{q^m}-[\ell]_{q^m})^k=q^{mk\ell}$.
Hadamard multiplication therefore replaces $t$ by $tq^{mk}$.
The last $i+1$ denominator factors cancel, leaving $(t;q)_{mk}$.
\end{proof}

\subsection{Product formulas and specializations}\label{sec:special}
In the following colored length specializations,
$\alpha=\alpha_{\E}$ and $b=b_{\E}$. The same arguments
apply to $\mathcal U$ with its own effective weight.

\subsubsection{Ranks divisible by the order of the root}
\begin{corollary}\label{cor:multiple}
With the notation of Theorem~\ref{thm:main},
\begin{equation}\label{eq:multiple}
 \E_{mk}^{(r)}(t,q,\zeta;\uu)
 =\left(\prod_{\substack{1\le j\le mk\\m\nmid j}}(1-tq^j)\right)
 \calB_k(t,q^m;b).
\end{equation}
\end{corollary}
\begin{proof}
For $i=0$, the first factor in~\eqref{eq:main} is $1/(1-t)$,
the identity for Hadamard multiplication. Furthermore,
\[
 \frac{(t;q)_{mk+1}}{(t;q^m)_{k+1}}
 =\prod_{\substack{1\le j\le mk\\m\nmid j}}(1-tq^j).
\]
\end{proof}

For rank $mk+1$, since $\B_1(t,q,\zeta;\alpha)=1+\alpha\zeta t$,
Theorem~\ref{thm:main} also gives
\begin{equation}\label{eq:remainderone}
 \frac{\E_{mk+1}^{(r)}(t,q,\zeta;\uu)}{(t;q)_{mk+2}}
 =\sum_{\ell\ge0}
 \bigl([\ell+1]_q+\alpha\zeta[\ell]_q\bigr)
 \bigl([\ell+1]_{q^m}+b[\ell]_{q^m}\bigr)^k t^\ell.
\end{equation}
Thus the first nonzero remainder requires no additional enumerator.

\subsubsection{The degenerate cycle}
\begin{corollary}\label{cor:degenerate}
If $(-\alpha)^m=1$, equivalently $b=-1$, then
\begin{equation}\label{eq:degenerate}
 \E_{mk+i}^{(r)}(t,q,\zeta;\uu)
 =(t;q)_{mk}\,
 \E_i^{(r)}(tq^{mk},q,\zeta;\uu).
\end{equation}
\end{corollary}
\begin{proof}
The factor raised to $k$ in~\eqref{eq:maincoeff} becomes
$q^{mk\ell}$. Hadamard multiplication by
$\sum_{\ell\ge0}q^{mk\ell}t^\ell$ replaces $t$ by $tq^{mk}$.
Hence the normalized enumerator is
\[
 \frac{\E_i^{(r)}(tq^{mk},q,\zeta;\uu)}
 {(tq^{mk};q)_{i+1}}.
\]
Multiplying by $(t;q)_{mk+i+1}$ proves the assertion.
\end{proof}

\subsubsection{Unweighted colors and their residue classes}
Write $\E_n^{(r)}(t,q,p)$ for~\eqref{eq:Edef} with every
$u_c=1$. Then $\alpha=[r-1]_\zeta$.

\begin{corollary}\label{cor:residue}
Fix a primitive $m$th root $\zeta$. If $r,s\ge1$ satisfy
$r\equiv s\pmod m$, then, for every $n\ge0$,
\begin{equation}\label{eq:residue}
 \E_n^{(r)}(t,q,\zeta)=\E_n^{(s)}(t,q,\zeta).
\end{equation}
If $r\equiv1\pmod m$, this common enumerator is
$A_n(t,q,\zeta)$. If $m\mid r$, then
\begin{equation}\label{eq:divisiblecolors}
 \E_{mk+i}^{(r)}(t,q,\zeta)
 =(t;q)_{mk}\E_i^{(r)}(tq^{mk},q,\zeta).
\end{equation}
\end{corollary}
\begin{proof}
Complete sums of $m$ consecutive powers of $\zeta$ vanish. Thus
$[r-1]_\zeta=[s-1]_\zeta$, and Lemma~\ref{lem:fold} proves
\eqref{eq:residue}. For $r\equiv1$, this value is zero, leaving
only the symmetric-group part of $\B_n$.
For $m\mid r$, it is $-\zeta^{-1}$, so
$(-\alpha)^m=1$ and Corollary~\ref{cor:degenerate} applies.
\end{proof}

The dependence only on $r\bmod m$ concerns the unweighted color
specialization. Independent color variables retain the more detailed
dependence through~\eqref{eq:effective}.

\subsubsection{Type B and the signed case}
For $r=2$ and unweighted colors, $\alpha=1$ and
$b=(-1)^{m+1}$. Consequently, for odd $m$,
\begin{equation}\label{eq:Bodd}
 \B_{mk}(t,q,\zeta;1)
 =\left(\prod_{\substack{1\le j\le mk\\m\nmid j}}(1-tq^j)\right)
 \calB_k(t,q^m;1),
\end{equation}
whereas for even $m$, all remainders satisfy
\begin{equation}\label{eq:Beven}
 \B_{mk+i}(t,q,\zeta;1)
 =(t;q)_{mk}\B_i(tq^{mk},q,\zeta;1).
\end{equation}
For $m=2$, the remainders are $i=0,1$ and
$\B_1(t,q,-1;1)=1-t$. Thus, for every $n\ge0$,
\begin{equation}\label{eq:Bsign}
 \sum_{\beta\in B_n}(-1)^{\ell_B(\beta)}
 t^{\des_B(\beta)}q^{\maj(\beta)}=(t;q)_n.
\end{equation}
This is the identity of Eu, Lin and Lo
\cite[Theorem 5.1]{ELL}, recovered here as the order-two case.
Their theorem also identifies this polynomial with the enumerator
weighted by $(-1)^{\negat(\beta)}$. Their Theorem 5.2 is recovered from the other statistic:
\begin{equation}\label{eq:Isign}
 I_n(t,q,-1;1)=
 \bigl(1+(-1)^{n-1}tq^{n-1}\bigr)(t;q)_{n-1},\qquad n\ge1.
\end{equation}
Indeed $I_1(t,q,p;1)=1+t$, whereas
$\B_1(t,q,p;1)=1+pt$, and Corollary~\ref{cor:BIdegenerate}
handles both parities. For example,
\[
 I_5(t,q,\mathrm i;1)=(t;q)_4(1+tq^4),\qquad
 \B_5(t,q,\mathrm i;1)=(t;q)_4(1+\mathrm i\,tq^4).
\]
At ranks divisible by the root order, Theorem~\ref{thm:B}
makes the two weighted enumerators equal; their remainder factors
can differ as in this example.

For unweighted colored groups at $p=-1$, the residue reduction gives
the following dichotomy:
\begin{equation}\label{eq:paritycolors}
 \E_n^{(r)}(t,q,-1)=
 \begin{cases}
  A_n(t,q,-1),&r\text{ odd},\\
  (t;q)_n,&r\text{ even}.
 \end{cases}
\end{equation}

\subsubsection{Two small higher-order examples}
\begin{example}[A cubic root in type B]\label{ex:cubic}
Let $\zeta^2+\zeta+1=0$. For $m=3,r=2$, the effective weight
is $b=1$. Since $\calB_1(t,Q;1)=1+t$,
\[
 \B_3(t,q,\zeta;1)=(1-tq)(1-tq^2)(1+t).
\]
For the same cubic root and $r=3$, one has
$\alpha=1+\zeta=-\zeta^2$, hence $b=-1$ and
\[
 \E_3^{(3)}(t,q,\zeta)=(t;q)_3.
\]
These illustrate the distinction between the ordinary type-B factor
and the degenerate cycle.
\end{example}

\begin{example}[Three colors at a fourth root]\label{ex:quartic}
Take $m=4$, $\zeta=\mathrm i$, and $r=3$, with unweighted colors.
Then $\alpha=1+\mathrm i$ and
$b=-(1+\mathrm i)^4=4$. Therefore
\begin{equation}\label{eq:quartic}
 \E_4^{(3)}(t,q,\mathrm i)
 =(1-tq)(1-tq^2)(1-tq^3)(1+4t).
\end{equation}
More generally, the large-rank factor in~\eqref{eq:main} is
$\calB_k(t,q^4;4)$. By color reduction at $p=1$, this is the
unweighted five-color Euler--Mahonian polynomial
$\E_k^{(5)}(t,q^4,1)$. Thus a cyclotomic specialization with three
colors can produce an ordinary enumerator with five colors.
\end{example}

\subsubsection{The Eulerian specialization}
Setting $q=1$ in~\eqref{eq:maincoeff} yields
\begin{equation}\label{eq:qone}
 [t^\ell]\frac{\E_{mk+i}^{(r)}(t,1,\zeta;\uu)}{(1-t)^{mk+i+1}}
 =\left([t^\ell]\frac{\E_i^{(r)}(t,1,\zeta;\uu)}{(1-t)^{i+1}}\right)
 \bigl(1+(1+b)\ell\bigr)^k.
\end{equation}
In particular,
\[
 \E_{mk}^{(r)}(t,1,\zeta;\uu)
 =(1-t)^{(m-1)k}\calB_k(t,1;b).
\]
These are polynomial identities and remain valid at special values
of the color variables where the effective parameter degenerates.

\subsection{Flag-major specialization and congruences}\label{sec:flag}

\subsubsection{The Biagioli--Zeng flag-major statistic}
Retain the order~\eqref{eq:order} and define, as in \cite{BZ},
\[
 \fmaj_G(\gamma)=r\maj(\gamma)+\col(\gamma).
\]
For a further color-weight variable $a$, put
\[
 F_n^{(r)}(t,x,p;a)=\sum_{\gamma\in G(r,n)}
 t^{\des_G(\gamma)}x^{\fmaj_G(\gamma)}p^{\ell_G(\gamma)}
 a^{\col(\gamma)}.
\]
It is obtained from~\eqref{eq:Edef} by $q=x^r$ and $u_c=(ax)^c$.
Theorem~\ref{thm:main} therefore gives the following formula.

\begin{corollary}\label{cor:flag}
With $n=mk+i$ and $\zeta$ as before, set
\[
 \alpha=ax[r-1]_{ax\zeta},\qquad b=-(-\alpha)^m.
\]
Then
\begin{equation}\label{eq:flag}
 \frac{F_{mk+i}^{(r)}(t,x,\zeta;a)}{(t;x^r)_{mk+i+1}}
 =\frac{F_i^{(r)}(t,x,\zeta;a)}{(t;x^r)_{i+1}}
 \otimes\frac{\calB_k(t,x^{rm};b)}{(t;x^{rm})_{k+1}}.
\end{equation}
\end{corollary}

The designation $\fmaj_G$ here fixes the convention. The flag-major
statistic defined from the color-block order used in
\cite[Section 2.3]{ELL} is pointwise different. Formula~\eqref{eq:flag}
concerns the statistic displayed above, together with $\des_G$;
it does not assert a formula for a flag-descent statistic based on
another order.

\subsubsection{Polynomial congruences}
The main identity can also be stated without choosing a root in
$\mathbb C$. Let $\Phi_m(p)$ be the $m$th cyclotomic polynomial and
put
\[
 \alpha(p)=\sum_{c=1}^{r-1}u_cp^{c-1},\qquad
 b_m(p)=-(-\alpha(p))^m.
\]
For each $n,\ell$, the coefficients
\[
 C_{n,\ell}(q,p;\uu)
 =[t^\ell]\frac{\E_n^{(r)}(t,q,p;\uu)}{(t;q)_{n+1}}
\]
are polynomials in $q,p,\uu$. By~\eqref{eq:maincoeff},
\begin{equation}\label{eq:congruence}
 C_{mk+i,\ell}(q,p;\uu)
 \equiv C_{i,\ell}(q,p;\uu)
 \bigl([\ell+1]_{q^m}+b_m(p)[\ell]_{q^m}\bigr)^k
 \pmod{\Phi_m(p)}.
\end{equation}
Indeed the difference vanishes at a primitive $m$th root, and
$\Phi_m$ is its minimal polynomial. Equivalently, the proof of
Theorem~\ref{thm:B} may be carried out in
$\mathbb Z[p,q,\uu]/(\Phi_m(p))$. This formulation is useful for
exact computation and does not require numerical approximations
to roots of unity.

\section{Fixed-point insertion and the joint kernel}\label{sec:insertion}

Insertion preserves length parity, but not length modulo higher
orders. To retain the exact change we attach a height to each gap.
Throughout this section $\delta=[\sigma_1^{c_1},\ldots,
\sigma_k^{c_k}]\in G(r,k)$ is a nonempty colored core, so
$k\ge1$; it need not be a derangement. For the empty core,
insertion gives only identity windows and both kernels are
$1/(1-x)$.

\subsection{Gap coordinates and crossing heights}
For $g=(g_0,\ldots,g_k)\in\mathbb Z_{\ge0}^{k+1}$, set
$N=\sum_hg_h$, $n=k+N$, and
\[
 T=\{t_1<\cdots<t_k\},\qquad
 t_a=a+\sum_{h=0}^{a-1}g_h.
\]
Construct $I_g(\delta)\in G(r,n)$ by putting
$t_{\sigma_a}^{c_a}$ in position $t_a$, and putting the uncolored
letter $j$ in every position $j\notin T$. These are the newly
inserted positive fixed points; any positive fixed points of the
core remain fixed. If the core is a derangement, this is a
bijection between gap vectors of sum $N$ and its derangement-part
fiber in rank $k+N$.
These gaps are the multiplicity form of the fixed-point insertion
slots in \cite[Section 3.1]{FZ}.

For $0\le h\le k$, define
\begin{equation}\label{eq:eta}
 \eta_h(\delta)=\#\{a\le h:\sigma_a>h\},
\end{equation}
and
\begin{equation}\label{eq:rho}
 \rho_h(\delta)
 =\#\{a\le h:\sigma_a>h\}
 +\#\{a>h:\sigma_a>h,\ c_a>0\}.
\end{equation}
For an ordinary derangement, the second term vanishes, and
$\rho_h$ is the usual crossing height of a permutation diagram.
The corresponding local fixed-point weight $p^{2\rho_h}$ appears
in the Laguerre-history model of \cite[Corollary 11]{SZ}.
The first term counts arrows crossing the cut after $h$ from left
to right. The second counts nonzero-colored entries whose position
and absolute value are both to the right of that cut. In particular,
$\rho_0(\delta)=\#\{a:c_a>0\}$ and $\rho_k(\delta)=0$.
For $X=X_L$ or $X_I$, use the height
\[
 H_h^X=\begin{cases}\rho_h,&X=X_L,\\\eta_h,&X=X_I.\end{cases}
\]
We also abbreviate these as $H_h^L$ and $H_h^I$.
In particular $H_k^X=0$.

\begin{theorem}[Exact insertion defect]\label{thm:defect}
For every $g\in\mathbb Z_{\ge0}^{k+1}$,
\begin{equation}\label{eq:defect}
 \ell_G(I_g(\delta))
 =\ell_G(\delta)+2\sum_{h=0}^kg_h\rho_h(\delta).
\end{equation}
For underlying unsigned inversion,
\begin{equation}\label{eq:defectI}
 \inv(|I_g(\delta)|)=\inv(|\delta|)+2\sum_{h=0}^kg_h\eta_h(\delta).
\end{equation}
Thus $X(I_g(\delta))-X(\delta)=2\sum_hg_hH_h^X$ for either $X$.
\end{theorem}
\begin{proof}
Comparisons between the $k$ core entries are preserved by the
increasing relabeling $a\mapsto t_a$. The color terms $c_a-1$
in the length formula are also preserved. Attribute each remaining
contribution to one inserted fixed point $f$ in gap $h$.

For an uncolored entry in position $t_a$ and with absolute value
$t_{\sigma_a}$, its inversion contribution with $f$ is
\[
 \bigl|\mathbf1_{a>h}-\mathbf1_{\sigma_a>h}\bigr|.
\]
For a nonzero-colored entry, that inversion contribution is
$\mathbf1_{a>h}$. In addition, increasing its absolute value from
$\sigma_a$ to $t_{\sigma_a}$ adds
$\mathbf1_{\sigma_a>h}$ to the length for this fixed point.
Thus coloring an entry replaces the displayed absolute difference
by the sum of its two indicators. The increase over the uncolored
contribution is $2\mathbf1_{a>h,\,\sigma_a>h}$.

For an ordinary permutation the number of arrows crossing this cut
in either direction is twice $\#\{a\le h:\sigma_a>h\}$, since
the numbers crossing left to right and right to left agree.
Consequently this fixed point contributes exactly
$2\rho_h(\delta)$. Fixed points in the same or different gaps
are increasing and create no mutual inversions. Summing their
contributions proves~\eqref{eq:defect}. For unsigned inversions
the two crossing directions have cardinality $\eta_h$, yielding
\eqref{eq:defectI}. The argument also allows fixed entries in
the core.
\end{proof}

\begin{corollary}[The inversion insertion series]\label{cor:lengthfiber}
For either inversion statistic $X$ and any colored core,
\begin{equation}\label{eq:lengthfiber}
 \sum_{g\in\mathbb Z_{\ge0}^{k+1}}x^{|g|}
 p^{X(I_g(\delta))-X(\delta)}
 =\prod_{h=0}^k\frac1{1-xp^{2H_h^X(\delta)}}.
\end{equation}
For a derangement core this is its derangement-part insertion
series. Let $\zeta$ have order $m\ge2$, put
$d=m/\gcd(m,2)$, $\omega=\zeta^2$, and let
$M_s=\#\{h:H_h^X(\delta)\equiv s\pmod d\}$.
At $p=\zeta$, the right side is
\begin{equation}\label{eq:heightresidues}
 \prod_{s=0}^{d-1}(1-x\omega^s)^{-M_s}.
\end{equation}
Its $x^N$ coefficient is a quasipolynomial in $N$ of period
dividing $d$ and degree at most $\max_sM_s-1$.
\end{corollary}
\begin{proof}
Sum the geometric series for the independent gaps. At a root
only $H_h^X\bmod d$ matters. Partial fractions express each
coefficient as a sum of terms
$\omega^{sN}\binom{N+j-1}{j-1}$, proving the last assertion.
\end{proof}

For $p=-1$, every gap weight is $1$, and the insertion series is
$(1-x)^{-k-1}$. Thus length parity is constant throughout each
fiber, recovering the parity fact used by Ji--Zhang. For higher
orders, the residue profile in~\eqref{eq:heightresidues} is the
additional datum required for the length enumeration.

\begin{example}\label{ex:defect21}
For the ordinary derangement $\delta=21$,
$\ell_G(\delta)=1$ and $(\rho_0,\rho_1,\rho_2)=(0,1,0)$.
Hence
\[
 \sum_{N\ge0}x^N
 \sum_{\operatorname{dp}(\pi)=21}p^{\inv(\pi)}
 =\frac{p}{(1-x)^2(1-xp^2)},
\]
where the inner sum is over rank $N+2$. At $p=\mathrm i$, this
becomes $\mathrm i/((1-x)(1-x^2))$, giving
\[
 \sum_{\substack{\pi\in\Sn_{N+2}\\\operatorname{dp}(\pi)=21}}
 \mathrm i^{\inv(\pi)}
 =\mathrm i\left(\left\lfloor\frac N2\right\rfloor+1\right).
\]
The middle gap accounts exactly for the extra factor $p^2$
identified in~\eqref{eq:fiberobstruction}.
\end{example}

\subsection{Descents and major index under insertion}
Let $d_h=\mathbf1_{h\in\Des_G(\delta)}$ for $0\le h<k$.
Set $A_0=0$ and, for $1\le h<k$, set
\[
 A_h=\mathbf1_{c_h=0,\,\sigma_h>h}.
\]
For $0\le h<k$, put
\[
 B_h=\mathbf1_{c_{h+1}>0\ \text{or}\ \sigma_{h+1}\le h}.
\]
In particular, $B_0=d_0$. Define
\[
 D_h(g)=\begin{cases}
 d_h,&g_h=0,\\
 A_h+B_h,&g_h>0.
 \end{cases}
\]
The two formulas agree at $h=0$.

\begin{lemma}\label{lem:gapdesc}
For $\pi=I_g(\delta)$,
\begin{align}
 \des_G(\pi)&=\sum_{h=0}^{k-1}D_h(g),\label{eq:gapdes}\\
 \maj(\pi)&=\sum_{h=0}^{k-1}
 \left(h+\sum_{j<h}g_j\right)D_h(g)
 +\sum_{h=0}^{k-1}B_hg_h.\label{eq:gapmaj}
\end{align}
\end{lemma}
\begin{proof}
An empty internal gap preserves its original comparison. A nonempty
gap replaces that comparison by the two boundaries of an increasing
block of fixed points. At the left boundary a descent occurs
precisely when the preceding entry is an uncolored excedance,
giving $A_h$. At the right boundary a descent occurs precisely
when the following entry is colored or its absolute value lies to
the left of the cut, giving $B_h$.

These boundaries occur in positions
$h+\sum_{j<h}g_j$ and $h+\sum_{j\le h}g_j$ respectively.
For the initial gap, the first boundary is the comparison with
zero and the same formulas apply. Appending fixed points after
the last entry creates no descent. Summing the descent counts
and their positions proves the result.
\end{proof}

For $S\subseteq\{1,\ldots,k-1\}$, define
\begin{align*}
 D_h(S)&=\begin{cases}A_h+B_h,&h\in S,\\d_h,&h\notin S,\end{cases}
 &D(S)&=\sum_{h=0}^{k-1}D_h(S),\\
 J(S)&=\sum_{h=0}^{k-1}hD_h(S),
 &L_h(S)&=B_h+\sum_{j=h+1}^{k-1}D_j(S).
\end{align*}
The subscript $\delta$ will be supplied when two derangements are
being considered.

\begin{theorem}[The Euler--Mahonian insertion kernel]\label{thm:kernel}
For an arbitrary colored core $\delta$ and $X=X_L$ or $X_I$, define
\begin{equation}\label{eq:kerneldef}
 K_\delta^X(x;t,q,p)=\sum_{g\in\mathbb Z_{\ge0}^{k+1}}
 x^{|g|}t^{\des_G(I_g(\delta))}q^{\maj(I_g(\delta))}
 p^{X(I_g(\delta))-X(\delta)}.
\end{equation}
Then
\begin{equation}\label{eq:kernel}
 K_\delta^X(x;t,q,p)=\frac1{1-x}
 \sum_{S\subseteq\{1,\ldots,k-1\}}
 \frac{t^{D(S)}q^{J(S)}}{1-xp^{2H_0^X}q^{D(S)}}
 \prod_{h\in S}\frac{xp^{2H_h^X}q^{L_h(S)}}
                       {1-xp^{2H_h^X}q^{L_h(S)}}.
\end{equation}
For a derangement core the gap sum is equivalently the sum over
its derangement-part fiber. We abbreviate $K_\delta=K_\delta^L$.
\end{theorem}
\begin{proof}
Fix the set $S$ of nonempty internal gaps. The other internal
gaps are zero, whereas $g_0,g_k\ge0$ are unrestricted. Equations
\eqref{eq:gapdes}--\eqref{eq:gapmaj} give
\[
 \des_G(\pi)=D(S),\qquad
 \maj(\pi)=J(S)+g_0D(S)+\sum_{h\in S}g_hL_h(S).
\]
Use~\eqref{eq:defect} or~\eqref{eq:defectI}, according to $X$,
and sum the independent gap variables.
The final gap gives $(1-x)^{-1}$ because $H_k^X=0$.
The initial gap gives the other denominator outside the product.
Each $g_h\ge1$ for $h\in S$ gives a geometric series starting
with its first power, producing the displayed numerator and
denominator. Summing over $S$ proves~\eqref{eq:kernel}.
\end{proof}

This is a rational formula for each individual fiber at arbitrary
$p$, and hence at every root of unity. Setting $t=q=1$ reduces
it to~\eqref{eq:lengthfiber}. It uses the cut heights and the local
excedance/subcedance data of $\delta$, rather than its major index
and length alone.

For $\delta=21$, the two supports give the explicit formula
\begin{equation}\label{eq:kernel21}
 pK_{21}(x;t,q,p)
 =\frac{ptq}{(1-x)(1-xq)}
 +\frac{p^3t^2q^3x}{(1-x)(1-xq^2)(1-xp^2q)}.
\end{equation}
Its $x$ coefficient is
$ptq+ptq^2+p^3t^2q^3$, accounting for $213,132,321$ with
their full $(\des,\maj,\inv)$ weights.

\subsection{Color-block flag-major kernels}
Let $K_{\delta,\prec}^X$ denote the same gap series as
$K_\delta^X$, using $\des_\prec,\maj_\prec$ in place of
$\des_G,\maj$. Compute its local data with
$d_h=\mathbf1_{h\in\Des_\prec(\delta)}$; the indicators
$A_h,B_h$ are unchanged, because all inserted letters have color
zero. The proof of Lemma~\ref{lem:gapdesc} therefore gives
\eqref{eq:kernel} with this replacement of $d_h$.
Define
\begin{equation}\label{eq:flagkernel}
 \mathcal H_\delta^X(x;t,q,p)
 =q^{\col(\delta)}K_{\delta,\prec}^X(x;t,q^r,p).
\end{equation}
It is the rational gap series retaining
$t^{\des_\prec(I_g(\delta))}q^{\fmaj_\prec(I_g(\delta))}$
and the relative $X$-weight. Formula~\eqref{eq:flagkernel}
applies to arbitrary cores, including those with positive fixed
points. For the empty core it gives $(1-x)^{-1}$.

\subsection{An exact number of fixed points}
\begin{corollary}\label{cor:exactfkernel}
Let $0\le f\le n$, and let $X=X_L$ or $X_I$. Then
\begin{align}
 &\sum_{\substack{\gamma\in G(r,n)\\\fix_0(\gamma)=f}}
 p^{X(\gamma)}w(\gamma)t^{\des_G(\gamma)}q^{\maj(\gamma)}
 \notag\\
 &\qquad=[x^f]\sum_{\delta\in\mathcal D_{n-f}^{(r)}}
 p^{X(\delta)}w(\delta)K_\delta^X(x;t,q,p).
 \label{eq:exactfkernel}
\end{align}
For color-block flag-major statistics, the same identity holds
with $\des_G,\maj,K$ replaced by
$\des_\prec,\fmaj_\prec,\mathcal H$.
\end{corollary}
\begin{proof}
Deleting all positive fixed points leaves a unique derangement
of rank $n-f$, and its gap vector has sum $f$. Conversely, inserting
such a vector into a derangement creates precisely $f$ positive
fixed points. Colors and their weights are unchanged.
\end{proof}

\subsection{Marking all fixed points by arbitrary-core insertion}
Put
\[
 F_n^X(s;t,q,p;\uu)=\sum_{\gamma\in G(r,n)}
 s^{\fix_0(\gamma)}p^{X(\gamma)}w(\gamma)
 t^{\des_G(\gamma)}q^{\maj(\gamma)}.
\]
\begin{theorem}[Fixed-point marking]\label{thm:markedkernel}
For $X=X_L$ or $X_I$,
\begin{equation}\label{eq:markedkernel}
 \sum_{n\ge0}F_n^X(s;t,q,p;\uu)x^n
 =\sum_{k\ge0}x^k\sum_{\delta\in G(r,k)}
 p^{X(\delta)}w(\delta)K_\delta^X((s-1)x;t,q,p).
\end{equation}
For color-block flag-major statistics, replace $K$ by
$\mathcal H$ and the corresponding statistics on the left.
\end{theorem}
\begin{proof}
For a colored window $\gamma$, choose any subset
$A\subseteq\Fix_0(\gamma)$. Deleting just $A$ gives a unique
colored core $\delta$ and gap vector $g$, with
$\gamma=I_g(\delta)$ and $|g|=|A|$. The core may still have
positive fixed points, which is why all of $G(r,k)$ occurs.
Conversely every core and gap vector determines such a pair.
The substitution $(s-1)x$ supplies $(s-1)^{|A|}$, and summing
over the subsets of $\Fix_0(\gamma)$ gives
$\sum_A(s-1)^{|A|}=s^{\fix_0(\gamma)}$. All other weights are
those of the resulting window. Each coefficient of $x$ is a
finite sum, so the formal rearrangement is valid.
\end{proof}

\section{Derangements and signed fixed-point identities}\label{sec:derangements}

\subsection{The derangement insertion expansion}
Let
\[
 D_n^X(t,q,p;\uu)=\sum_{\gamma\in\mathcal D_n^{(r)}}
 p^{X(\gamma)}w(\gamma)t^{\des_G(\gamma)}q^{\maj(\gamma)}.
\]
Setting $s=0$ in Theorem~\ref{thm:markedkernel} gives
\begin{equation}\label{eq:derangementGF}
 \sum_{n\ge0}D_n^X(t,q,p;\uu)x^n
 =\sum_{k\ge0}x^k\sum_{\delta\in G(r,k)}
 p^{X(\delta)}w(\delta)K_\delta^X(-x;t,q,p).
\end{equation}
This is inclusion-exclusion with explicit rational insertion
weights. It applies at every $p$, and the flag-major version is
obtained by replacing $K$ by $\mathcal H$.
At $t=q=1$ it simplifies to
\begin{equation}\label{eq:lengthderangementGF}
 \sum_{n\ge0}D_n^X(1,1,p;\uu)x^n
 =\sum_{k\ge0}x^k\sum_{\delta\in G(r,k)}
 \frac{p^{X(\delta)}w(\delta)}
 {\prod_{h=0}^k(1+xp^{2H_h^X(\delta)})}.
\end{equation}
For type B and $n\ge1$, the negative-entry-refined signed
balances are
\begin{align}
 \sum_{\beta\in\mathcal D_n^B}(-1)^{\ell_B(\beta)}z^{\negat(\beta)}
 &=(-1)^n\bigl(1+n(z-1)\bigr),\label{eq:signedcountL}\\
 \sum_{\beta\in\mathcal D_n^B}(-1)^{\inv(|\beta|)}z^{\negat(\beta)}
 &=(-1)^{n-1}\bigl(n(1+z)-1\bigr).\label{eq:signedcountI}
\end{align}
Indeed the full-group signed sums vanish in ranks at least two.
Their rank-one values are $1-z$ and $1+z$, respectively, and
\eqref{eq:lengthderangementGF} at $p=-1$ reduces to
$1/(1+x)+x(1\mp z)/(1+x)^2$.

\subsection{Conventions and the signed refinement}
Ji and Zhang \cite{JZ} use the type-B flag-major index
$\fmaj_\prec=2\maj_\prec+\negat$, where
\begin{equation}\label{eq:Jiorder}
 -1\prec-2\prec\cdots\prec-n\prec0\prec1\prec\cdots\prec n.
\end{equation}
This differs from the natural-order major index used here, and
from the Biagioli--Zeng flag-major index in Section~\ref{sec:flag}.
For example, on the five signed derangements
$\mathcal D_2^B=\mathcal D_2^{(2)}$ of rank two,
\begin{equation}\label{eq:threeindices}
 \sum_{\beta\in\mathcal D_2^B}(-1)^{\ell_B(\beta)}q^{\maj(\beta)}=q,
 \qquad
 \sum_{\beta\in\mathcal D_2^B}(-1)^{\ell_B(\beta)}q^{\fmaj_\prec(\beta)}
 =q+q^3-q^4.
\end{equation}
Their type-D formulas also use intrinsic type-D length and
major-index conventions.

\subsection{A fixed-point-compatible color reduction}
Color reduction preserves the positive fixed-point set. For independent
variables $v_0,\ldots,v_{n-1}$ and $x_1,\ldots,x_n$, its proof gives
\begin{align}
 &\sum_{\gamma\in G(r,n)}p^{\ell_G(\gamma)}
 \prod_{d\in\Des_G(\gamma)}v_d
 \prod_{j=1}^nu_{c_j}\prod_{j\in\Fix_0(\gamma)}x_j\notag\\
 &\qquad=\sum_{\beta\in B_n}p^{\ell_B(\beta)}
 \alpha(p)^{\negat(\beta)}
 \prod_{d\in\Des_B(\beta)}v_d
 \prod_{j\in\Fix_+(\beta)}x_j.
 \label{eq:markedfold}
\end{align}
Indeed folding changes neither the positive fixed points nor the
descent set. Consequently the residue-class identity
\eqref{eq:residue} holds after derangement restriction, and also
with arbitrary positive fixed-point-set markers. This follows from
color reduction alone; it makes no assertion about factorization in
the rank.

\begin{theorem}[Fixed-point-set refinement of Eu--Lin--Lo]\label{thm:fixedset}
For every $n\ge0$,
\begin{align}
 &\sum_{\beta\in B_n}(-1)^{\ell_B(\beta)}
 \prod_{d\in\Des_B(\beta)}v_d
 \prod_{j\in\Fix_+(\beta)}x_j\notag\\
 &\quad=\sum_{\beta\in B_n}(-1)^{\negat(\beta)}
 \prod_{d\in\Des_B(\beta)}v_d
 \prod_{j\in\Fix_+(\beta)}x_j
 =\prod_{j=1}^n(x_j-v_{j-1}).\label{eq:fixedsetproduct}
\end{align}
\end{theorem}
\begin{proof}
We retain the positive fixed-point set in the involution used in
\cite[proof of Theorem 5.1]{ELL}. If $|\beta_j|=j$ for every $j$,
leave the window unchanged. Otherwise let $h$ be the least index
with $|\beta_h|\ne h$, and let $s$ satisfy $|\beta_s|=h$.
Then $s>h$. Change the sign of $\beta_s$.

This operation is an involution because it leaves the absolute-value
window unchanged. It also preserves $\Fix_+$: the modified letter
has absolute value $h\ne s$. All letters at positions at least $h$
other than $s$ have absolute value greater than $h$. Hence the
neighbors of position $s$, when present, have absolute value greater than $h$,
and replacing $h$ by $-h$ or conversely changes neither neighboring
comparison. Thus the entire descent set is preserved. The negative
entry count changes by one, and the length changes by
$\pm(2h-1)$, as in \cite{ELL}. Both signs in
\eqref{eq:fixedsetproduct} are therefore reversed on nonfixed pairs.

The surviving windows satisfy $|\beta_j|=j$. For each $j$, choosing
$+j$ gives the factor $x_j$. Choosing $-j$ gives a descent at $j-1$
and the factor $-v_{j-1}$. Their independent choices give the product.
\end{proof}

\begin{corollary}\label{cor:fixnumber}
If $\alpha(-1)=\sum_{c=1}^{r-1}(-1)^{c-1}u_c=1$, then
\begin{equation}\label{eq:coloredfixedset}
 \sum_{\gamma\in G(r,n)}(-1)^{\ell_G(\gamma)}
 \prod_{d\in\Des_G(\gamma)}v_d
 \prod_{j=1}^nu_{c_j}\prod_{j\in\Fix_0(\gamma)}x_j
 =\prod_{j=1}^n(x_j-v_{j-1}).
\end{equation}
In particular, for unweighted colors and even $r$,
\begin{equation}\label{eq:fixnumber}
 \sum_{\gamma\in G(r,n)}(-1)^{\ell_G(\gamma)}
 t^{\des_G(\gamma)}q^{\maj(\gamma)}s^{\fix_0(\gamma)}
 =\prod_{j=0}^{n-1}(s-tq^j),
\end{equation}
and for $0\le f\le n$,
\begin{equation}\label{eq:exactfixed}
 \sum_{\substack{\gamma\in G(r,n)\\\fix_0(\gamma)=f}}
 (-1)^{\ell_G(\gamma)}t^{\des_G(\gamma)}q^{\maj(\gamma)}
 =(-t)^{n-f}q^{\binom{n-f}{2}}\qbinom{n}{f}_q.
\end{equation}
In particular,
\begin{equation}\label{eq:signedderangements}
 \sum_{\gamma\in\mathcal D_n^{(r)}}(-1)^{\ell_G(\gamma)}
 t^{\des_G(\gamma)}q^{\maj(\gamma)}
 =(-t)^nq^{\binom n2}\qquad(r\text{ even}).
\end{equation}
\end{corollary}
\begin{proof}
Apply~\eqref{eq:markedfold} and Theorem~\ref{thm:fixedset}.
For even $r$ and $u_c=1$, the alternating sum $\alpha(-1)$ is $1$.
Set $v_j=tq^j$ and $x_j=s$, and use the finite $q$-binomial theorem
to extract $s^f$. Setting $s=0$ gives the derangement identity.
\end{proof}

These formulas retain descent information absent from the univariate
Mahonian derangement enumerators. They are consequences of the
Eu--Lin--Lo involution and the color reduction; their proof does not
use a cyclotomic deletion formula.

At $t=q=1$, the even-color case recovers the signed derangement
balance $(-1)^n$. For odd $r$, folding gives $\alpha(-1)=0$
and hence the ordinary signed balance $(-1)^{n-1}(n-1)$ for $n\ge1$.
These numerical balances agree with the type-A and type-B cases
in~\cite{JZ} and the colored counts in~\cite[Theorem 3.9]{AAA},
whose length convention is different. We include them only as
specializations of the fixed-point refinement.

\subsection{Signed flag-major formulas with an exact fixed-point number}
Define, with the color-block convention~\eqref{eq:blockflag},
\begin{align*}
 d_{r,k}^-(q;\uu)&=\sum_{\delta\in\mathcal D_k^{(r)}}
 (-1)^{\ell_G(\delta)}q^{\fmaj_\prec(\delta)}w(\delta),\\
 M_{r,j}^-(q;\uu)&=\sum_{\gamma\in G(r,j)}
 (-1)^{\ell_G(\gamma)}q^{\fmaj_\prec(\gamma)}w(\gamma).
\end{align*}
Both rank-zero polynomials are one.

\begin{theorem}\label{thm:flagfixed}
For $0\le f\le n$ and every $r\ge1$,
\begin{equation}\label{eq:flagfixed}
 \sum_{\substack{\gamma\in G(r,n)\\\fix_0(\gamma)=f}}
 (-1)^{\ell_G(\gamma)}q^{\fmaj_\prec(\gamma)}w(\gamma)
 =\qbinom{n}{f}_{q^r}d_{r,n-f}^-(q;\uu).
\end{equation}
In particular,
\begin{equation}\label{eq:flagderangement}
 d_{r,k}^-(q;\uu)=\sum_{a=0}^k
 (-1)^a q^{r\binom a2}\qbinom{k}{a}_{q^r}
 M_{r,k-a}^-(q;\uu).
\end{equation}
\end{theorem}
\begin{proof}
The colored major-index insertion identity \cite[Section 3]{FZ}
and \cite[Theorem 3.1 and its general-color extension]{CGH} gives
\[
 \sum_{\substack{\gamma\in G(r,n)\\\operatorname{dp}(\gamma)=\delta}}
 q^{\fmaj_\prec(\gamma)}
 =q^{\fmaj_\prec(\delta)}\qbinom{n}{k}_{q^r},
 \qquad\delta\in\mathcal D_k^{(r)}.
\]
Theorem~\ref{thm:defect} shows that length parity is constant
on this fiber, for every $r$. Its color weight is also constant.
Multiplying by these constant factors and summing over rank
$n-f$ derangements proves~\eqref{eq:flagfixed}. Summing over
$f$ and applying $q$-binomial inversion with base $q^r$ gives
\eqref{eq:flagderangement}.
\end{proof}

Thus $r=1$ recovers the type-A exact-fixed-point refinement of
Ji--Zhang's Proposition 2.6, and $r=2$ recovers its type-B
counterpart from Proposition 3.6. The unweighted unsigned relation
is already part of the fix-Euler--Mahonian framework
\cite{FZ,CGH}; here parity preservation supplies the signed
length weight.

\begin{proposition}[The signed full-group flag-major factor]
\label{prop:flagfull}
Set
\[
 C_r(q;\uu)=\sum_{c=0}^{r-1}u_c^2q^{2c},\qquad
 A_r(q;\uu)=\sum_{c=0}^{r-1}(-1)^cu_cq^c.
\]
For $j\ge0$,
\begin{equation}\label{eq:flagfull}
 M_{r,j}^-(q;\uu)=
 C_r(q;\uu)^{\lfloor j/2\rfloor}
 A_r(q;\uu)^{j\bmod2}
 \prod_{a=1}^j[a]_{(-1)^{a-1}q^r}.
\end{equation}
For unweighted colors, $C_r=[r]_{q^2}$ and $A_r=[r]_{-q}$.
\end{proposition}
\begin{proof}
Introduce the auxiliary color-block polynomial
\[
 P_j(v,Q,p;\mathbf a)=\sum_{\gamma\in G(r,j)}
 v^{\des_F(\gamma)}Q^{\maj_\prec(\gamma)}
 p^{\inv(|\gamma|)}\prod_{h=1}^j a_{c_h},
\]
where $\des_F$ counts only descent positions $1,\ldots,j-1$.
Ordinary compatible decreasing sequences have one block for
each level $0\le h\le\ell$ and color $0\le c<r$.
Every such block has increasing absolute labels in the color-block
order. The argument of Proposition~\ref{prop:coeffI} therefore gives
\begin{equation}\label{eq:flagbuckets}
 [v^\ell]\frac{P_j(v,Q,p;\mathbf a)}{(v;Q)_{j+1}}
 =\sum_{\sum_{h,c}m_{h,c}=j}
 \qbinom{j}{(m_{h,c})}_p
 \prod_{h,c}(a_cQ^h)^{m_{h,c}}.
\end{equation}
Write $j=2b+i$, $i\in\{0,1\}$. At $p=-1$, the no-carry
evaluation of Lemma~\ref{lem:Lucas} makes this coefficient
\[
 \left(\sum_ca_c^2\right)^b\left(\sum_ca_c\right)^i
 [\ell+1]_{Q^2}^{,b}[\ell+1]_Q^{,i}.
\]
The last two factors are the same normalized coefficient for the
ordinary signed Euler--Mahonian polynomial $A_j(v,Q,-1)$.
Consequently
\[
 P_j(v,Q,-1;\mathbf a)=
 \left(\sum_ca_c^2\right)^b\left(\sum_ca_c\right)^i
 A_j(v,Q,-1).
\]
Now take $v=1$, $Q=q^r$ and $a_c=(-1)^cu_cq^c$.
The parity identity~\eqref{eq:coloredparity} identifies the left
side with $M_{r,j}^-$. The Gessel--Simion signed major-index
product, as recalled in \cite{AGR,JZ}, gives
$A_j(1,Q,-1)=\prod_{a=1}^j[a]_{(-1)^{a-1}Q}$.
This proves~\eqref{eq:flagfull}.
\end{proof}

\begin{corollary}[Unweighted natural-order parity reduction]
\label{cor:naturalparityfixed}
Let $k=n-f$. For the Biagioli--Zeng major index,
\begin{equation}\label{eq:naturalparityfixed}
 \sum_{\substack{\gamma\in G(r,n)\\\fix_0(\gamma)=f}}
 (-1)^{\ell_G(\gamma)}q^{\maj(\gamma)}
 =\qbinom{n}{f}_q
 \begin{cases}
 d_{1,k}^-(q),&r\text{ odd},\\
 (-1)^kq^{\binom k2},&r\text{ even}.
 \end{cases}
\end{equation}
The even-color formula retains descent number by the additional
factor $t^k$, as in~\eqref{eq:exactfixed}.
\end{corollary}
\begin{proof}
For odd $r$, color folding at $p=-1$ gives
$\alpha(-1)=[r-1]_{-1}=0$. Since it preserves the positive
fixed-point set, only ordinary permutations remain. Their
major-index fiber formula is Ji--Zhang's Proposition 2.6.
For even $r$, use~\eqref{eq:exactfixed} with $t=1$.
\end{proof}

At $q=1$, both conventions give the signed fixed-point count
\[
 \binom nf(-1)^{n-f}
 \begin{cases}1-(n-f),&r\text{ odd},\\1,&r\text{ even}.
 \end{cases}
\]
This includes $f=n$, with value one.

\subsection{Parity of descent positions and alternating descents}
Pan, Xu and Zeng \cite{PXZ} use the order~\eqref{eq:order} and
length~\eqref{eq:length} to enumerate odd and even descent positions.
Write
\[
 \odes_G(\gamma)=\#\{j\in\Des_G(\gamma):j\text{ odd}\},\qquad
 \edes_G(\gamma)=\#\{j\in\Des_G(\gamma):j\text{ even}\}.
\]
Position zero is included among the even positions. In their notation
$\col_G(\gamma)$ counts nonzero-colored entries, whereas our
$\col(\gamma)$ is their $\operatorname{csum}(\gamma)$.
To keep these conventions distinct, put
$\nu(\gamma)=\#\{j:c_j>0\}$ and denote their polynomial, with $p$
as the length variable, by
\[
 P_n^{(r)}(x,y;a,p)=\sum_{\gamma\in G(r,n)}
 x^{\odes_G(\gamma)}y^{\edes_G(\gamma)}
 a^{\nu(\gamma)}p^{\ell_G(\gamma)}.
\]
Thus their color weight corresponds to $u_c=a$ for every $c>0$,
giving $\alpha=a[r-1]_\zeta$ in~\eqref{eq:effective}.
Moreover, $\des_G=\odes_G+\edes_G$ and
$\maj\equiv\odes_G\pmod2$, so, writing
$a\mathbf1=(a,\ldots,a)$ for the color vector,
\begin{align}
 \E_n^{(r)}(t,1,p;a\mathbf1)&=P_n^{(r)}(t,t;a,p),\notag\\
 \E_n^{(r)}(t,-1,p;a\mathbf1)&=P_n^{(r)}(-t,t;a,p).
 \label{eq:paritybridge}
\end{align}
These specializations relate the two refinements. General independent
odd/even descent weights retain information beyond these evaluations.

\begin{corollary}[Parity and fixed-point refinement]\label{cor:parityfixed}
Let $r$ be even and $n\ge0$. Then
\begin{align}
 &\sum_{\gamma\in G(r,n)}(-1)^{\ell_G(\gamma)}
 x^{\odes_G(\gamma)}y^{\edes_G(\gamma)}
 q^{\maj(\gamma)}s^{\fix_0(\gamma)}\notag\\
 &\qquad=
 \prod_{j=0}^{\lceil n/2\rceil-1}(s-yq^{2j})
 \prod_{j=0}^{\lfloor n/2\rfloor-1}(s-xq^{2j+1}).
 \label{eq:parityfixed}
\end{align}
Define $\operatorname{altdes}_G$ to count odd descents and even
ascents, including a possible ascent at zero. Then
\begin{equation}\label{eq:alternatingfixed}
 \sum_{\gamma\in G(r,n)}(-1)^{\ell_G(\gamma)}
 z^{\operatorname{altdes}_G(\gamma)}s^{\fix_0(\gamma)}
 =(s-z)^{\lfloor n/2\rfloor}(sz-1)^{\lceil n/2\rceil}.
\end{equation}
\end{corollary}
\begin{proof}
In~\eqref{eq:coloredfixedset}, take every $u_c=1$, every fixed-point
marker equal to $s$, and $v_j=yq^j$ for even $j$, $v_j=xq^j$ for
odd $j$. Evenness of $r$ ensures $\alpha(-1)=1$, proving
\eqref{eq:parityfixed}. Since there are $\lceil n/2\rceil$ even
positions in $\{0,\ldots,n-1\}$,
\[
 \operatorname{altdes}_G
 =\odes_G+\lceil n/2\rceil-\edes_G.
\]
Set $q=1$, $x=z$, $y=z^{-1}$ in~\eqref{eq:parityfixed} and
multiply by $z^{\lceil n/2\rceil}$. This proves
\eqref{eq:alternatingfixed} as a polynomial identity.
\end{proof}

At $s=1$, equation~\eqref{eq:alternatingfixed} is
\[
 \sum_{\gamma\in G(r,n)}(-1)^{\ell_G(\gamma)}
 z^{\operatorname{altdes}_G(\gamma)}
 =(-1)^{\lceil n/2\rceil}(1-z)^n,
\]
which is precisely \cite[Theorem 2.2(1)]{PXZ}. Color reduction and
the Eu--Lin--Lo involution therefore give a combinatorial proof of
this even-color case, as sought in \cite[Remark (iii) after
Theorem 2.2]{PXZ}. At $s=0$, the same refinement gives
\[
 \sum_{\gamma\in\mathcal D_n^{(r)}}(-1)^{\ell_G(\gamma)}
 z^{\operatorname{altdes}_G(\gamma)}
 =(-1)^n z^{\lfloor n/2\rfloor}.
\]

\subsection{Why fixed-point deletion does not extend directly}
The signed fiber identities in
\cite[Propositions 2.6 and 3.6]{JZ} use preservation of length
parity under this deletion, together with a shuffle identity. In
type A the identity is
\begin{equation}\label{eq:JZfiber}
 \sum_{\substack{\pi\in\Sn_n\\\operatorname{dp}(\pi)=\sigma}}
 (-1)^{\inv(\pi)}q^{\maj(\pi)}
 =(-1)^{\inv(\sigma)}q^{\maj(\sigma)}\qbinom{n}{k}_q,
 \qquad\sigma\in\mathcal D_k^{(1)}.
\end{equation}
The corresponding type-B identity uses $\fmaj_\prec$ and the
Gaussian coefficient with base $q^2$.

For a primitive root of order $m>2$, the literal analogue of
\eqref{eq:JZfiber} already fails on $\Sn_3$. The fiber above
$\sigma=21$ consists of $213,132,321$. Its weight is
\[
 \zeta q+\zeta q^2+\zeta^3q^3,
\]
whereas replacing the sign on the right-hand side by
$\zeta^{\inv(\sigma)}$ gives $\zeta q(1+q+q^2)$.
The difference is
\begin{equation}\label{eq:fiberobstruction}
 \zeta(\zeta^2-1)q^3\ne0.
\end{equation}
Deleting the fixed middle entry of $321$ changes the inversion
number from $3$ to $1$. This preserves parity, but not its residue
modulo $m>2$.

Nor does Theorem~\ref{thm:main} hold after replacing every
enumerator by its derangement restriction. For $r=m=n=2$ and
$\zeta=-1$, formula~\eqref{eq:signedderangements} gives the actual
polynomial $t^2q$. A direct derangement substitution in
Corollary~\ref{cor:multiple} would instead yield
$-t(1-tq)$, since the rank-one weighted type-B derangement
enumerator at weight $b=-1$ is $-t$. This counterexample separates
the full-group cyclotomic theorem from its signed derangement
consequences. The height-weighted expansion~\eqref{eq:derangementGF} retains
this information. Compressing it to a factor depending only on
rank remains a separate question.

For Ji--Zhang's actual flag-major convention, the smallest
nonempty fiber gives the stronger type-B example
\begin{equation}\label{eq:minusonefiber}
 \sum_{N\ge0}x^N\sum_{\substack{\pi\in B_{N+1}\\
                         \operatorname{dp}(\pi)=(-1)}}
 p^{\ell_B(\pi)}q^{\fmaj_\prec(\pi)}
 =\frac{pq}{(1-x)(1-xp^2q^2)}.
\end{equation}
It follows directly from $H_0^L=1$, $H_1^L=0$ and
\eqref{eq:flagkernel}. At $p=\zeta$, its rank-two coefficient
is $\zeta q+\zeta^3q^3$, differing from a literal Gaussian
replacement by $\zeta(\zeta^2-1)q^3$.

The full-group congruence~\eqref{eq:rankcongruence} also fails
under derangement restriction. In rank two, the natural-order
signed length polynomial is $t^2q$, while the unsigned-inversion
sign polynomial is $t^2q-2t(1+q)$. Thus the two full-group
factorizations and the insertion formulas have distinct scopes.

\section{Color restrictions and further questions}\label{sec:restrictions}

\subsection{Position-dependent color choices in the signed case}
Allow a separate color weight $u_{j,c}$ in position $j$, with
$u_{j,0}=1$. The proof of~\eqref{eq:markedfold} still applies,
but the type-B weight in position $j$ becomes
\[
 \alpha_j(p)=\sum_{c=1}^{r-1}u_{j,c}p^{c-1}.
\]
If $\alpha_j(-1)=1$ for every $j$, the folded signed
enumerator has no residual position-dependent negative-entry
weight. Consequently~\eqref{eq:coloredfixedset}, including the
entire descent set and positive fixed-point set, remains valid.
In particular, position $j$ may allow the unweighted colors
$0,1,\ldots,2s_j-1$, with independently chosen $s_j\ge1$
and $2s_j\le r$. This gives a class of position restrictions
for which the signed product survives.

\subsection{Filtering by total color}
Let $d\mid r$ and use the inherited $G(r,1,n)$ statistics on
\[
 G(r,d,n)=\{\gamma\in G(r,n):\col(\gamma)\equiv0\pmod d\}.
\]
Denote its enumerator by $E_n^{(r,d)}(t,q,p)$. Let $\eta$ be a
primitive $d$th root of unity. The elementary root filter gives
\begin{equation}\label{eq:subgroupfilter}
 E_n^{(r,d)}(t,q,p)
 =\frac1d\sum_{j=0}^{d-1}
 E_n^{(r)}(t,q,p;a=\eta^j),
\end{equation}
where the scalar notation means $u_c=a^c$.
For a primitive $m$th root $\zeta$, put
\[
 \alpha_j=\eta^j[r-1]_{\eta^j\zeta},\qquad
 b_j=-(-\alpha_j)^m.
\]
Applying Theorem~\ref{thm:main} inside the filter gives
\begin{align}
 \frac{E_{mk+i}^{(r,d)}(t,q,\zeta)}{(t;q)_{mk+i+1}}
 =\frac1d\sum_{j=0}^{d-1}
 \left(\frac{E_i^{(r)}(t,q,\zeta;a=\eta^j)}{(t;q)_{i+1}}
 \otimes\frac{\calB_k(t,q^m;b_j)}{(t;q^m)_{k+1}}\right).
 \label{eq:subgroupHadamard}
\end{align}
In particular,
\begin{equation}\label{eq:subgroupmultiple}
 E_{mk}^{(r,d)}(t,q,\zeta)
 =\left(\prod_{\substack{1\le h\le mk\\m\nmid h}}(1-tq^h)\right)
 \frac1d\sum_{j=0}^{d-1}\calB_k(t,q^m;b_j).
\end{equation}
These formulas provide the color-sum-subgroup extension for
inherited statistics. When $G(2,2,n)=D_n$, they do not replace
the intrinsic type-D descent and major-index conventions used
in Ji--Zhang's formulas.

\subsection{Insertion in color-sum subgroups}
Let $\mathcal D_k^{(r,d)}=\mathcal D_k^{(r)}\cap G(r,d,k)$.
The inherited statistics admit the same insertion formulas as the
ambient colored group.

\begin{corollary}[Subgroup insertion fibers]\label{cor:subgroupfiber}
For $d\mid r$ and $\delta\in\mathcal D_k^{(r,d)}$, every gap
insertion $I_g(\delta)$ lies in $G(r,d,k+N)$, where $N=\sum_hg_h$.
The one-sided and two-sided insertion kernels restricted to these
subgroups are therefore exactly~\eqref{eq:kernel}
and~\eqref{eq:inversekernel}. In particular,
\begin{equation}\label{eq:subgroupfiber}
 \sum_{N\ge0}x^N
 \sum_{\substack{\pi\in G(r,d,k+N)\\\operatorname{dp}(\pi)=\delta}}
 p^{\ell_G(\pi)-\ell_G(\delta)}
 =\prod_{h=0}^{k}\frac1{1-xp^{2\rho_h(\delta)}}.
\end{equation}
For the empty derangement the series is $1/(1-x)$.
\end{corollary}
\begin{proof}
Inserted fixed points have color zero and the other colors are
unchanged. Hence $\col(I_g(\delta))=\col(\delta)$ as integers.
The congruence defining $G(r,d,n)$ is therefore preserved, and the
entire ambient insertion fiber lies in the subgroup. Apply
Theorem~\ref{thm:kernel} and Corollaries~\ref{cor:lengthfiber}
and~\ref{cor:inversekernel}.
\end{proof}

\subsection{Further questions}

The proof of Theorem~\ref{thm:main} removes carries by multinomial
evaluation, and each complete cycle supplies $1-(-\alpha)^m$.
A direct combinatorial description for $m>2$ would explain the
effective weight on the type-B factor.

The remaining insertion problem is to sum the kernels of
Theorem~\ref{thm:kernel} over derangements without retaining every
cut profile. At a root of order $m$, the required length data are
the heights modulo $m/\gcd(m,2)$; local descent data must also be
controlled. Intrinsic type-D statistics additionally require the
last-entry refinement in Ji--Zhang's treatment of
$\Delta_n=\{\beta\in B_n:\beta_n>0\}$.

Beyond the fiber refinement in Corollary~\ref{cor:inversekernel}
and the position restrictions and subgroup filters in
Section~\ref{sec:restrictions}, an arbitrary-order rank
factorization with general position weights or simultaneous direct
and inverse Eulerian statistics remains to be established.

The fake-degree framework of Barcelo--Reiner--Stanton~\cite{BRS}
suggests a uniform extension beyond the family $G(r,d,n)$.
Such an extension would require an Eulerian refinement compatible
with the chosen length and major-index statistics, and a substitute
for insertion coordinates in exceptional reflection groups.

\appendix

\section{Inverse statistics use the same heights}\label{app:inverse}
Here $K_\delta$ denotes the length kernel $K_\delta^L$.

Positive fixed points are shared by a colored permutation and its
group inverse. In gap coordinates,
\begin{equation}\label{eq:inversegaps}
 I_g(\delta)^{-1}=I_g(\delta^{-1}),\qquad
 \rho_h(\delta^{-1})=\rho_h(\delta).
\end{equation}
The first identity follows from using the same increasing
relabeling of positions and values. For the second, the two
crossing counts of the underlying permutation are equal, and the
nonzero-colored entries with both coordinates to the right of
a cut are unchanged by transposing the coordinates. Inversion
replaces a nonzero color by its additive inverse modulo $r$,
which remains nonzero.

\begin{corollary}[Two-sided insertion kernel]\label{cor:inversekernel}
Let $K_\delta^{\leftrightarrow}(x;t,q,u,v,p)$ be the insertion
series with the additional weight
$u^{\des_G(\pi^{-1})}v^{\maj(\pi^{-1})}$. Use superscripts
$+$ and $-$ for the quantities $D,J,L$ computed from $\delta$
and $\delta^{-1}$ respectively. Then
\begin{align}
 K_\delta^{\leftrightarrow}
 &=\frac1{1-x}\sum_{S\subseteq\{1,\ldots,k-1\}}
 \frac{t^{D^+(S)}u^{D^-(S)}q^{J^+(S)}v^{J^-(S)}}
 {1-xp^{2\rho_0}q^{D^+(S)}v^{D^-(S)}}\notag\\
 &\hspace{2em}\times\prod_{h\in S}
 \frac{xp^{2\rho_h}q^{L_h^+(S)}v^{L_h^-(S)}}
 {1-xp^{2\rho_h}q^{L_h^+(S)}v^{L_h^-(S)}}.
 \label{eq:inversekernel}
\end{align}
\end{corollary}
\begin{proof}
Apply Lemma~\ref{lem:gapdesc} to $\delta$ and $\delta^{-1}$
with the same gap vector, then repeat the geometric summation
in Theorem~\ref{thm:kernel}.
\end{proof}

Thus inverse descent and major-index refinements are available
at the fiber level. A closed rank-reduction theorem for the sum
over all derangements requires further compression of these data.

For the unsigned-inversion kernel, the same formulas hold after
replacing $\rho_h$ by $\eta_h$. The equality of the two crossing
directions proves $\eta_h(\delta^{-1})=\eta_h(\delta)$, so the
same geometric summation applies.

\section{Galois-twisted color-block flag-major kernels}\label{app:galois}
For an arbitrary finite complex reflection group $W$,
Barcelo, Reiner and Stanton~\cite{BRS} define the bimahonian
polynomial
\[
 \mathcal M_W^\theta(U,V)
 =\sum_{\chi\in\operatorname{Irr}(W)}
 f^{\theta(\chi)}(U)f^{\bar\chi}(V),
\]
where $f^\chi$ is the fake-degree polynomial, $\bar\chi$ is the
contragredient character, and $\theta$ is a Galois automorphism
of a cyclotomic field of definition. Their Theorems 1.3 and 1.4
give a Galois-twisted flag-major interpretation for wreath products
and root-of-unity evaluations via regular elements, respectively.
Their flag-major statistic uses the color-block order
\eqref{eq:blockorder} and index~\eqref{eq:blockflag}.
Biagioli and Caselli~\cite{BC} extend descent and flag-major
enumeration to projective reflection groups $G(r,p,s,n)$.
Here we retain the inherited color-block statistics on $G(r,d,n)$
and add the Biagioli--Zeng length on each insertion fiber.
Throughout this subsection $\ell_G$ and $\rho_h$ retain their
definitions~\eqref{eq:length} and~\eqref{eq:rho}; in particular,
length is still computed using~\eqref{eq:order}.

For an integer $a$ with $\gcd(a,r)=1$, let $\theta_a$ multiply
every color by $a$ modulo $r$, taking residues in $\{0,\ldots,r-1\}$.
This is the entrywise Galois automorphism of the monomial matrices.
It preserves $G(r,d,n)$, since
$\col(\theta_a(\gamma))\equiv a\col(\gamma)\pmod d$.

\begin{proposition}[Galois-twisted flag-major insertion kernel]
\label{prop:Galflag}
Let $d\mid r$, $\gcd(a,r)=1$, and
$\delta\in\mathcal D_k^{(r,d)}$ with $k\ge1$. Set
$\tau=\theta_a(\delta^{-1})$ and $\rho_h=\rho_h(\delta)$.
Compute the local quantities $A_h,B_h,D,J,L$ of
Section~\ref{sec:insertion} for $\delta$ and $\tau$, replacing
only $d_h$ by $\mathbf1_{h\in\Des_\prec}$.
Use superscripts $+$ and $-$ for these quantities on $\delta$
and $\tau$, respectively. Define
\begin{align*}
 \mathcal H_{\delta,a}(x;t,U,u,V,p)
 &=\sum_{N\ge0}x^N
 \sum_{\substack{\pi\in G(r,d,k+N)\\\operatorname{dp}(\pi)=\delta}}
 t^{\des_\prec(\pi)}u^{\des_\prec(\theta_a(\pi^{-1}))}\\
 &\hspace{2em}\times
 U^{\fmaj_\prec(\pi)}V^{\fmaj_\prec(\theta_a(\pi^{-1}))}
 p^{\ell_G(\pi)-\ell_G(\delta)}.
\end{align*}
Then
\begin{align}
 \mathcal H_{\delta,a}
 &=\frac{U^{\col(\delta)}V^{\col(\tau)}}{1-x}
 \sum_{S\subseteq\{1,\ldots,k-1\}}
 \frac{t^{D^+(S)}u^{D^-(S)}U^{rJ^+(S)}V^{rJ^-(S)}}
 {1-xp^{2\rho_0}U^{rD^+(S)}V^{rD^-(S)}}\notag\\
 &\hspace{2em}\times\prod_{h\in S}
 \frac{xp^{2\rho_h}U^{rL_h^+(S)}V^{rL_h^-(S)}}
 {1-xp^{2\rho_h}U^{rL_h^+(S)}V^{rL_h^-(S)}}.
 \label{eq:Galflag}
\end{align}
The empty-derangement kernel is $1/(1-x)$.
\end{proposition}
\begin{proof}
Increasing relabeling preserves comparisons in the color-block
order, so an empty internal gap retains its descent indicator.
For a nonempty gap, the fixed points are uncolored. A left boundary
descent occurs exactly when the preceding letter is an uncolored
excedance, and a right boundary descent occurs exactly when the
following letter is colored or lies to the left of the cut.
Thus $A_h,B_h$ are unchanged, and Lemma~\ref{lem:gapdesc} holds
with the indicated replacement of $d_h$. The initial and final
gaps obey the same rules as before.

Inversion and the color automorphism commute with insertion:
\[
 \theta_a(I_g(\delta)^{-1})=I_g(\tau),\qquad
 \rho_h(\tau)=\rho_h(\delta).
\]
The second equality uses~\eqref{eq:inversegaps} and the fact that
$\theta_a$ preserves zero and nonzero colors. Moreover,
$\col(I_g(\delta))=\col(\delta)$ and
$\col(I_g(\tau))=\col(\tau)$ as integers. Therefore replacing
the two major-index variables by $U^r,V^r$ and multiplying by
$U^{\col(\delta)}V^{\col(\tau)}$ converts the geometric summation
in Corollary~\ref{cor:inversekernel} to~\eqref{eq:Galflag}.
Corollary~\ref{cor:subgroupfiber} justifies the subgroup restriction.
\end{proof}

For $d=1$, setting $t=u=p=1$ and summing the coefficients
$[x^{n-k}]\mathcal H_{\delta,a}$ over all ranks $k\le n$ and all
$\delta\in\mathcal D_k^{(r)}$ recovers
$\mathcal M_{G(r,1,n)}^{\theta_a}(V,U)$ from
\cite[Theorem 1.3]{BRS}. Formula~\eqref{eq:Galflag} additionally
retains length on each fiber. It does not provide a closed
rank-reduction formula for the sum over all derangements.

\end{document}